\PassOptionsToPackage{dvipsnames}{xcolor}
 \documentclass[a4paper,12pt]{amsart}
\usepackage[latin1]{inputenc}
\usepackage[T1]{fontenc}
\usepackage{amsfonts,amssymb}
\usepackage[normalem]{ulem}
\usepackage{float} 
\usepackage[pagebackref]{hyperref} 
\usepackage{pstricks}
\usepackage{amsthm}
\usepackage[all]{xy}
\usepackage[pdftex]{graphicx}
\usepackage{graphicx}  
 \usepackage[first=0,last=9]{lcg}
 \usepackage{tikz}
 \usetikzlibrary{tikzmark}
 \usepackage{ytableau}
\usetikzlibrary{calc, positioning, fit, shapes.misc}

 \usepackage{color, colortbl}

 \tikzset{every picture/.style={remember picture}}
\definecolor{Gray}{gray}{0.9}
\newcommand{\n}[0]{{\mathbf N}}

\renewcommand{\S}[0]{{\mathcal S}}

\definecolor{Beaver}{rgb}{0.4,0.3,0.28}
\definecolor{Shadow}{rgb}{0.54,0.47,0.36}
\definecolor{grullo}{rgb}{0.66,0.6,0.53}

\newtheorem{introthm}{Theorem}

\newtheorem{theo}{Theorem}[section]

\newtheorem{prop}[theo]{Proposition}

\newtheorem{lem}[theo]{Lemma}
\newtheorem{conje}[theo]{Conjecture} 

\newtheorem{rem}[theo]{Remark}

\title[Partition identities associated with $A_r$-Surface singularities]{Partition identities associated with $A_r$ Surface singularities}
\author{Pooneh Afsharijoo, Pedro D. Gonz\'alez P\'erez and\\ Hussein Mourtada}

\thanks{\textbf{Acknowledgment.} The research of Pooneh Afsharijoo was funded by a \textit{Mar\'{\i}a Zambrano Postdoctoral fellowship}. 
{The first and second author were funded also by grants
PID2020-114750GB-C32 and PID2024-156181NB-C32 of MICIU/AEI
MCIN/AEI
/10.13039/501100011033 and 
FEDER, UE}}
\begin {document}

\maketitle

\begin{abstract} We prove a family of partition identities involving 
integer partitions in three colors. The conditions imposed on the 
types of partitions appearing in these identities involve 
constraints that arise in the Rogers-Ramanujan and Andrews-Gordon 
identities, as well as in their recent extensions. The identities 
established in this paper are associated with the $A_r$ surface 
singularities via the arc HP-series, which provides a measure of 
singularities of algebraic varieties defined using arc spaces.
\end{abstract}

\footnote{keywords: rational double point singularities; integer partitions; graded algebras; Hilbert series; differential algebra; arc spaces; cyclic quotient singularities; Gordon's identities; Rogers-Ramanujan identities}


\section{Introduction}

An (ordinary) partition $\lambda$ of a positive integer $n$ is a non-increasing sequence of positive integers $(\lambda_1 \geq \cdots \geq \lambda_{\ell})$ such that $\lambda_1+\cdots+\lambda_\ell = n$. The integers $\lambda_i$ are called the parts of $\lambda$. We denote by $L(\lambda)$ the number of parts (that is, the length) of $\lambda$, and by $p(n)$ the total number of partitions of $n$.

A fundamental family of partition identities for ordinary partitions, which plays a central role in this article, is given by Gordon's identities (see Theorem~$1$ in \cite{G}):

\begin{introthm} \label{Gordon}  (\textit{Gordon's identities}). Given integers $r\geq 2$ and $1\leq i \leq r,$ let $\mathcal{B}_{r,i}(n)$ denote the set of partitions of $n$ of the form $(\lambda_1,\dots, \lambda_s)$, where $\lambda_{j}-\lambda_{j+r-1} \geq 2$ and at most $i-1$ parts equal to $1$ and denote its cardinal by $B_{r,i}(n)$. Let $\mathcal{A}_{r,i}(n)$ denote the number of partitions of $n$ into parts $\not\equiv 0,\pm i \ ( \text{mod}.  2r+1)$ and denote its cardinal by $A_{r,i}(n)$. Then $A_{r,i}(n)=B_{r,i}(n)$ for all integers $n$.
\end{introthm}

Note that these identities generalize the so-called Rogers-Ramanujan identities, which correspond to the case $r=2$ in the above theorem. These identities have played an important role in several areas, including statistical mechanics, combinatorics and number theory, representation theory, probability theory, graph theory, as well as algebraic geometry and commutative algebra
\cite{AN,MM,AM2025,F,Bax,ADJM1,ADJM2,BMS,BIS,GIS,GOW,LZ,Bull}.

In this paper, we focus on integer partitions whose parts may occur in several colors (in fact, three). Let $c \geq 1$ be an integer. A partition is called \textit{$c$-colored} if each part can appear in $c$ different colors. Ordinary integer partitions correspond to the case $c=1$. For $c=3$, we consider the three colors black, red, and green. We denote by $i_b$ (respectively $i_r$ and $i_g$) the part of $i$ color black (respectively red and green).

The $3$-colored partitions of $2$ are:
$$2_b, 2_r, 2_g, 1_b+1_b, 1_b+1_r, 1_b+1_g, 1_r+1_g, 1_r+1_r, 1_g+1_g.$$
For a $3$-colored partition $\lambda$, we may consider the subpartitions $\lambda_b$, $\lambda_r$, and $\lambda_g$, obtained respectively from the black, red, and green parts of $\lambda$. For example, for the partition $\lambda = 7_b + 6_r + 3_b + 3_r + 1_r$, we have $\lambda_b = 7 + 3$, $\lambda_r = 6 + 3 + 1$, and $\lambda_g = \emptyset$.
\bigskip

To state the main theorem, we introduce new classes of integer partitions in $3$ colors, denoted by $\mathcal{F}_{r}$, indexed by an integer  $r \geq 2$.

\bigskip 
We first introduce some notation. For a $3$-colored partition $\lambda$, we denote by
\begin{itemize}
\item $\ell_b$ (respectively $\ell_r$ and $\ell_g$) the number of black (respectively red and green) parts of $\lambda$;
\item $k$ the smallest black part of $\lambda$ (if it exists); 
\item $i_k$ the $k$-th smallest red part of $\lambda$, when $\ell_r \geq k$.
\end{itemize}

\vspace{0.5cm}
We define $\mathcal{F}_{r}$ as the set of $3$-colored partitions $\lambda$ such that $\lambda_g \in \mathcal{B}_{r,r}(n)$ for some $n$. In addition, $\lambda$ satisfies one of the following conditions:
\begin{enumerate}
\item[1.] $\lambda = \emptyset$; or $\lambda$ has only one color; or $\lambda$ has exactly two colors, one of which is green.
\item[2.] $\lambda$ has exactly two colors, black and red, and
\begin{enumerate}
\item either $\ell_r \leq k-1$,
\item or $\ell_r \geq k$ and $k + i_k \geq r$.
\end{enumerate}
\item[3.] $\lambda$ has exactly three colors, and
\begin{enumerate}
\item either $\ell_r \leq k-1$,
\item or $\ell_r \geq k$, $k + i_k \geq r$, and $\ell_g - \#\{1_g\} < k + i_k - r + 1$,
\end{enumerate}
\end{enumerate}
where $\#\{1_g\}$ denotes the number of green parts which are equal to one.
\bigskip

We denote by $\n$ the set of nonnegative integers.
For $n \in \n$ and $r \geq 2$, let $F_r(n)$ denote the cardinality of $\mathcal{F}_r(n)$, the set of integer partitions of $n$ belonging to $\mathcal{F}_r$. The main result of this article is the following family of identities:
\begin{introthm} \label{th-frn}
For all $n \in \n$ and all integers $r \geq 2$, the numbers $F_r(n)$ are equal:
\[
F_2(n) = F_3(n) = F_4(n) = \cdots.
\]
Moreover, this common value equals the number of integer partitions of $n$ in which the part $1$ may appear in $3$ colors, while all other positive integers may appear in $2$ colors.
\end{introthm}

In Gordon's identities
each $r$ gives rise to a different identity, expressing an equality between the numbers of two types of partitions that satisfy conditions varying with $r$. What is remarkable in  Theorem  \ref{th-frn}
is that all the $r$-dependent conditions, appearing in the definition of $\mathcal{F}_r$ (and some of these conditions come from the defining conditions of $\mathcal{B}_{r,r}(n)$ in Gordon's identities)   are unified into a single identity containing a countably infinite number of terms. \\


This theorem originates from the study of certain surface singularities, namely the $A_r$ singularities, with $r \in \mathbf{N}$ and $r \geq 1$. We now explain how the $3$-colored partitions of type $\mathcal{F}_r$ are inspired by these surface singularities. Let $\mathbf{K}$ be a field of characteristic zero, and consider the polynomial 
\[
h_r:=z^{r}-xy\in \mathbf{K}[x,y,z].
\] 
 For $r \geq 2$, the surface defined by $X_{r-1}=\{h_r=0\}$ has, at the origin, what is called an $A_{r-1}$-singularity, which is a rational double point singularity; this singularity can be also seen as the quotient of $\mathbb{A}^2$  by the action of a cyclic group.\\\

The arc space $X_{r,\infty}$ of $X_r$ is the moduli space parameterizing the arcs on $X_r$. An arc $\gamma$ is given by a vector of power series

 \[\gamma(t)=(x(t),y(t),z(t))=(\sum_{i\in \mathbf{Z}_{\geq 0}}x_it^i,\sum_{i\in \mathbf{Z}_{\geq 0}}y_it^i,\sum_{i\in \mathbf{Z}_{\geq 0}}z_it^i)\in (\mathbf{K}[[t]])^3\] 
satisfying
$$h_r(\gamma(t))=0.$$
If we write the expansion 
\[h_r(\gamma(t))=\sum H_{i}t^i,
\]
then the arc space is the algebraic variety (actually a scheme) defined in the infinite-dimensional affine space with coordinates $x_i, y_i, z_i$ for $i \ge 0$, by the equations $$H_i=0, i\in \mathbf{Z}_{\geq 0}.$$ And the ring $\mathcal{A}_r$ of polynomial functions on $X_{r,\infty}$ is given by
$$\mathcal{A}_r=\frac{\mathbf{K}[x_i,y_i,z_i, i\in \mathbf{Z}_{\geq 0}]}{(H_i, i\in \mathbf{Z}_{\geq 0})}.$$

Up to a change of variables (replace respectively $x_i,y_i,z_i$ by $x_i/(i!),y_i/(i!),z_i/(i!)),$ we can show that the transform of the ideal $(H_i \mid i \in \mathbf{Z}_{\ge 0})$ is a differential ideal of $\mathbf{K}[x_i, y_i, z_i \mid i \in \mathbf{Z}_{\ge 0}]$.
More precisely, the ring $\mathbf{K}[x_i, y_i, z_i \mid i \in \mathbf{Z}_{\ge 0}]$, is naturally equipped with the derivation $D$ defined by \[
 D(x_i)=x_{i+1},~~~ D(y_i)=y_{i+1}~~~ D(z_i)=z_{i+1}. \]
The transform of the ideal $(H_i \mid i \in \mathbf{Z}_{\ge 0})$ by the change of variables satisfies that $D(f)\in (H_i \mid i \in \mathbf{Z}_{\ge 0})$ for any $f\in (H_i \mid i \in \mathbf{Z}_{\ge 0})$ (\textit{e.g.}, Proposition 2.3 in \cite{dif}); in the sequel we consider this ideal in the variables which are adapted to the differential structure; by abuse of notation, it will also be denoted by $(H_i \mid i \in \mathbf{Z}_{\ge 0})$.

Now, if we assign to each variable $x_i, y_i, z_i$ the weight $i$ for $i \in \mathbf{Z}_{\ge 0}$, then each polynomial $H_j$ is weighted homogeneous of weight $j$. This endows the ring $\mathcal{A}_r$ with a natural grading, \textit{i.e.}, 
$$\mathcal{A}_r=\bigoplus_{j\in \mathbf{Z}_{\geq 0}}\mathcal{A}_j,$$ 
where $\mathcal{A}_i \cdot\mathcal{A}_j\subset \mathcal{A}_{i+j}.$ We are particularly interested in the arcs $\gamma$ whose center is  singular point $0 \in X_{r-1},$ \textit{i.e.}, the space of arcs $\gamma$ satisfying $\gamma(0)=0;$ we denote it by $X_{r,\infty}^0$. 
The ring of functions on $X_{r,\infty}^0$ is the quotient 
$$\mathcal{A}^0: =
\mathbf{K}[x_i,y_i,z_i, i\in \mathbf{Z}_{\geq 1}]/
\mathfrak{a},
$$
where 
\begin{equation} \label{eq-ideal-a}
\mathfrak{a}=(\overline{H}_i, i\in \mathbf{Z}_{\geq 1})\subset \mathbf{K}[x_i,y_i,z_i,i \geq 1], 
\end{equation}
and $\overline{H}_i$ is obtained from $H_i$ by putting $x_0=y_0=z_0=0$.
\\


The ring $\mathcal{A}^0=\bigoplus_{j\in \mathbf{Z}_{\geq 0}}\mathcal{A}^0_j$ is again graded and one can consider its Hilbert series which is by definition the arc Hilbert Poincar\'e series (arc HP-series) of $X_{r-1}$ at $0:$
$$\mathrm{AHP}_{X_{r-1},0}(q):=\sum_{j\in \mathbf{Z}_{\geq 0}}
\mathrm{dim}_\mathbf{K}\mathcal{A}^0_jq^j.$$
Now, one observes that the algebra 
\[ 
\mathcal{S}:=\mathbf{K}[x_i,y_i,z_i, i\in \mathbf{Z}_{\geq 1}]
\] 
also inherits a grading from the weights assigned to the variables. To each monomial in $\mathcal{S}$, one can associate a $3$-colored partition by assigning to $x_i$, $y_i$, and $z_i$ a part of size $i$ colored black, red, and green, respectively. For example, the $3$-colored partition corresponding to the monomial $x_2^2y_2z_3z_4$ is $4_g+3_g+2_r+2_b+2_b.$\\

On the one hand, to determine the dimensions of the vector spaces $\mathcal{A}^0_j$, $j \in \mathbf{Z}_{\ge 0}$, we look for monomials in $\mathcal{S}$ whose images in $\mathcal{A}^0_j$ form the desired bases. A conjecture about which monomials serve this purpose arises from computations with Groebner base; this conjecture led us to introduce the families of integer partitions that we denote by $\mathcal{F}_r$, $r \in \mathbf{Z}_{\ge 2}$.
\\

On the other hand, the geometry of the jet schemes allows us to compute the arc-HP series 
 $\mathrm{AHP}_{X_{r-1},0}(q)$ \cite{M1}, which turns out to be the generating series for the number of integer partitions of $n$ with $3$ colors for the part $1$ and $2$ colors for all other positive integers. This observation leads directly to the statement of the main result of this paper, whose proof relies on several recent developments in combinatorics and differential algebra; in particular, \emph{the proof relies on a refinement of the main result of \cite{ADJM1}, which we establish using the results of \cite{BMS,PE} (see Theorem~\ref{refinement}).}




\section{A possible initial ideal of the space of arc at an $A_{r-1}$ singularity }

In this section, we are interested in  the description of  a Groebner basis of 
the ideal $\mathfrak{a} $ of the ring  $ \S$
(see  \eqref{eq-ideal-a}), with respect to the weighted reverse lexicographic order. More precisely, we focus on its initial (or leading) ideal. We will formulate below a conjecture describing this initial ideal. 

The ideals defined below are included in the ring  $\S$; 
let
 $$I_r:=(z_i^{r},z_i^{r-1}z_{i+1},\cdots,z_{i}z_{i+1}^{r-1}| \  i\geq 1)$$ 
 and denote by $J_r$ the ideal generated by following monomials:
 \begin{itemize}
\item[-] The monomials of $I_r;$
\item[-] The monomials of the form $$x_k y_{i_1}\cdots y_{i_k}$$ where $1\leq k, 1\leq i_1 \leq \cdots \leq i_k$ and $k+i_k\leq r-1;$
\item[-] The monomials of the form
$$x_k y_{i_1}\cdots y_{i_k} z_{j_1}\cdots z_{j_{i_k+k-r+1}},$$
where $1\leq k, 1\leq i_1 \leq \cdots \leq i_k, 2\leq j_1 \leq \cdots \leq j_{k+i_k-r+1}$ and $k+i_k\geq r.$
\end{itemize}

We have the following conjecture
\begin{conje}\label{conj} The initial ideal $\mathfrak{a}$ with respect to the weighted reverse lexicographical ordering is the ideal $J_r$ defined above.
\end{conje}

As explained in the introduction, a monomial in the graded ring  $\S$
 is naturally associated with a  three colored partition and vice versa. The partitions in $\mathcal{F}_r$ are exactly those associated with the monomials which are not in the ideal $J_r,$ or equivalently the monomials which are non-zero in the quotient ring $\S/J_r.$ In particular, the generating series of the sequence $F_r(n)$ is equal to the Hilbert series of $\S/J_r,$ because the monomials of $\S$ of weight $n$ which are not in $J_r$ form a basis of the weighted-homogeneous component of $\S/J_r$ of weight $n$. Now, it is a standard result in commutative algebra that the Hilbert series of an ideal coincides with the Hilbert series of its leading ideal with respect to a monomial ordering that respects the grading  (e.g., see \cite{GP}, the appendix of \cite{BMS}). Therefore, if Conjecture~\ref{conj} holds, it follows that the Hilbert series of $\S/J_r$ is equal to the Hilbert series of $\S/\mathfrak{a}$, which by definition is equal to $AHP_{X_{r-1}}(q).$ We know from \cite{M1} that

\begin{theo}\label{RDP}\[AHP_{X_{r-1}}(q)=\frac{1}{1-q^3}\prod_{i \geq 2}\frac{1}{1-q^i}\]

\end{theo}
  
One observes that the series in Theorem~\ref{RDP} is the generating series of integer partitions of $n$ in which the part $1$ may appear in $3$ colors, while all other positive integers may appear in $2$ colors. In particular, Conjecture~\ref{conj} implies Theorem{ \ref{th-frn}}.
However, in general it is very difficult to obtain a Groebner basis for an ideal that is not finitely generated, which is the case for our ideal $\mathfrak{a}.$ This explains the motivation behind Theorem{ \ref{th-frn}},
and the fact that it holds provides strong supporting evidence for the conjecture.


\bigskip 
Recall from the introduction that, for a $3$-colored partition $\lambda$, we denote by $\ell_b$ (respectively by $\ell_r$ and $\ell_g$) the number of black (respectively red and green) parts of $\lambda$. We 
 denote by $k$ the smallest black part of $\lambda$, if it exists.   If $\ell_r \geq k$, we denote by  $i_k$ the $k$-th smallest red part of $\lambda$. The proof of the main theorem suggests the following slightly different formulation (from that given in the introduction) of the defining conditions for the partitions in $\mathcal{F}_{r}$, which is the set of $3$-colored partitions $\lambda$ satisfying one of the following conditions:\\
\begin{enumerate}
\item[1.] Either $\lambda=\emptyset$ or all parts of $\lambda$ are colored in a single color: black, red or green; and if this color is green, then $\lambda \in \mathcal{B}_{r,r}(n)$ of Gordon's identities for some positive integer $n.$
\item[2.] $\lambda$ has two colors, either black and green or red and green, with green sub-partition $\lambda_g$ belonging to $\mathcal{B}_{r,r}(n)$ of Gordon's identities for some integer $n$.
\item[3.] $\lambda$ has two colors black and red, with 
\begin{enumerate}
\item either $\ell_r \leq k-1,$
\item or $\ell_r \geq k$ and $k+i_k\geq r.$ 
\end{enumerate}
\item[4.] $\lambda$ has three colors (black, red and green) with green sub-partition $\lambda_g$ belonging to $\mathcal{B}_{r,r}(n)$ of Gordon's identities for some integer $n$ and with
\begin{enumerate}

\item either $\ell_r\leq k-1,$
\item or $\ell_r \geq k$ and $ k+i_k\geq r$ with $\ell_g-\#\{1_g\}<k+i_k-r+1.$ 

\end{enumerate}
\end{enumerate}

Therefore, the generating series of $F_{r}(n)$ is equal to \emph{the sum of the generating series of these four types} of partitions.  Let us introduce first some notation about certain generating series.
\\

 We denote the generating series of all partions by
\[
\mathbf{H} := \sum_{n\in \mathbb{N}}p(n)q^n.
\] 
It satisfies the well-known identity:
$$\mathbf{H}=\prod_{j\geq 1} \frac{1}{1-q^j}.$$
We denote by $\mathbf{G}_r$ the generating series 
\[
\mathbf{G}_r := \sum_{n\in \mathbb{N}} B_{r,r}(n) q^n,
\]
where $B_{r,r}(n)$ is the number of partitions of $n$ appearing in Gordon's identities.
\\

The generating series of the numbers of partitions of type $1$ is obviously equal to:
\begin{equation}\label{S1}
S_1:=\underbrace{1}_{\text{empty set}}+\underbrace{2(\mathbf{H}-1)}_{\text{single-colored black or red partitions}}+\underbrace{(\mathbf{G}_r-1)}_{\text{single-colored green partitions}}
\end{equation}
 
It is also direct to compute the generating series of partitions of type $2$:

\begin{equation}\label{S2}
S_2:=2(\mathbf{H}-1)(\mathbf{G}_r-1).
\end{equation}
Determining the generating series for the number of partitions in $\mathcal{F}_r$ of the other types is more involved and requires a few formulas, which we will prove in the following sections.

\section{Useful generating series}

In this section, we provide the generating series for certain types of partitions that will be used later to prove our main theorem.

Let $k \geq 1$ be an integer. 
It is known that the generating function of non-empty partitions of length at most $k$ is (see \cite{A}):
\begin{equation}\label{red1}
P_{\leq k}:=\frac{1}{(q)_{k}}-1,
\end{equation}
 where $(q)_0:=1$ and $(q)_n: =(1-q)\cdots(1-q^n),$ for all integers $n\geq 1$.

To see a commutative algebraic proof of this formula, note that these partitions are those associated with monomials in the graded algebra $$\frac{\mathbf{K}[y_1,y_2,\cdots]}{(y_{i_1}\cdots y_{i_{k+1}})},$$
whose Hilbert-Poincar\'e series is equal to $\frac{1}{(q)_{k}}-1$ (see Lemma $3.2.$ in \cite{AM}).
\bigskip

We now have the following lemma, which provides the generating series for the black sub-partitions of $\mathcal{F}_{r}$ of type $3$ and type $4$:

\begin{lem}\label{black} Let $k\geq 1$ be an integer. The generating function of the partitions whose smallest part is equal to $k$ is:
$$\frac{q^k}{\prod_{j\geq k}(1-q^j)}.$$
\end{lem}
\begin{proof}
On the one hand, we can break down each  partition $\lambda=(\lambda_1\geq \cdots \geq \lambda_m \geq k)$  into two sub-partitions: $\mu=(\lambda_1\geq \cdots \geq \lambda_m)$ and $\eta=(k).$ On the one hand, since all parts of $\mu$ are at least equal to $k$, its generating function is $\frac{1}{\prod_{j\geq k}(1-q^j)}$. On the other hand $\eta$ is generated by $q^k.$
\end{proof}

\begin{rem} The previous lemma can be proved using commutative algebra. Note that the set of partitions $\mu$ in the proof of this lemma is in bijection with the set of monomials in $\mathbf{K}[x_k, x_{k+1}, \dots]$, whose Hilbert-Poincar\'e series is $\frac{1}{\prod_{j \geq k} (1-q^j)}$.

\begin{lem}\label{H} For the integer $m\geq 1$ we have:
$$\sum_{k \geq m} q^k(q)_{k-1}=(q)_{m-1} - \frac{1}{\mathbf{H}}{\red .} $$
\end{lem}

\end{rem}
\begin{proof}
The proof is by induction on $m.$ Recall that $\mathbf{H}$ is the generating series of all partitions. These are the empty partition and the partitions with smallest part equal to $k$ for some $k\geq 1.$ Therefore by Lemma \ref{black} we have:
$$\mathbf{H}=1+\sum_{k\geq 1}\frac{q^k}{\prod_{j\geq k}(1-q^j)}=1+\sum_{k\geq 1}\frac{q^k(q)_{k-1}}{\prod_{j\geq 1}(1-q^j)}=1+\mathbf{H}\sum_{k\geq 1} q^k(q)_{k-1},$$
Thus the initial case of the induction holds. Suppose now that the formula of the proposition is true for any integer $m'$ between $1$ and $m.$ For the case $m'=m+1$ we have:
$$(q)_m-\sum_{k\geq m+1} q^k(q)_{k-1}=(q)_m+q^m (q)_{m-1}-\sum_{k\geq m} q^k(q)_{k-1}=(q)_{m-1}-\sum_{k\geq m} q^k(q)_{k-1},$$
which, by induction hypothesis, is equal to $1/ \mathbf{H}.$

\end{proof}
The following lemma gives us the generating series of the red sub-partitions of $\mathcal{F}_{r}$ of type $3.b$ and type $4.b:$
\begin{lem}\label{red2} Let $k, i_k \geq 1$ be integers. The generating series of the partitions of length at least $k$ whose $k$-th smallest part is equal to $i_k$ is
$$\mathbf{H} \, q^{k+i_k-1} \frac{(q)_{k+i_k-2}}{(q)_{k-1}}.$$
\end{lem}
\begin{proof} On the one hand, we can decompose each partition $\lambda=(i_\ell \geq \cdots \geq i_1)$ with $\ell \geq k$ into three sub-partitions: $\mu=(i_\ell\geq \cdots \geq i_k), \gamma=(i_{k-1}-1 \geq \cdots \geq i_1-1)$ and $  \eta=(\underbrace{1,\cdots, 1}_{(k-1) \ \text{times}}).$  
\\
The partition $\eta$ is generated by $q^{k-1}$ and by Lemma \ref{black} the partition $\mu$ is generated by $$\frac{q^{i_k}}{\prod_{j\geq i_k}(1-q^j)}.$$
Note that $\gamma$ is a partition whose Young diagram fits into a rectangle of size $(i_k-1)\times (k-1)$ which is generated by the $q$-binomial coefficient: $$\left[{k+i_k-2\atop k-1}\right]_q := \frac{(q)_{k+i_k-2}}{(q)_{k-1}(q)_{i_k-1}}.$$
Multiplying these three generating series gives us the result.

\end{proof}
Let $\ell, n \geq 0$ and $r\geq 2$ be integers. Let $\mathcal{G}_{r,\ell}(n)$ be the set of partitions of $\mathcal{B}_{r,r}(n)$ whose number of non-equal to one parts is less than or equal to $\ell.$ Denote the cardinal of this set by $G_{r,\ell}(n).$ In the Proposition \ref{G} we give the generating series of these partitions.  Note that this is related to the generating series of the green sub-partitions of $\mathcal{F}_{r}$ of type $4.b$ (see Remark \ref{lamdag} below). Before representing this $q$-series, we will describe some key tools that we will use in order to prove it.\\

In \cite{Af}, the first author conjectured a new member of the Gordons identities. In other words, she defined a new set of partitions $\mathcal{C}_{r,i}(n)$ whose cardinal $C_{r,i}(n)$ is equal to the number $B_{r,i}(n)$ of Gordon's identities. Recently, this conjecture was proved in \cite{ADJM1} (see also \cite{ADJM2}, and \cite{PE} using a differential algebra approach). In order to prove this conjecture the authors of \cite{ADJM1} define another set of partitions $\mathcal{D}_{r,i}(n)$ and prove that it is equal to the set $\mathcal{C}_{r,i}(n)$ (see Theorem $2.3$ in \cite{ADJM1}). 
\\

In order to introduce $\mathcal{D}_{r,i}(n)$, we need to define the \textit{Durfee square} and the \textit{horizontal Durfee rectangle} of a partition:
The Durfee square (respectively the horizontal Durfee rectangle) of a partition is the largest square (respectively the largest rectangle of size $d \times (d+1)$) that can be placed in the top-left corner of its Young diagram.
Similarly, we can define successive Durfee squares/rectangles by drawing the first square/rectangle, then the next one below the first one and so on.\\

\begin{figure}
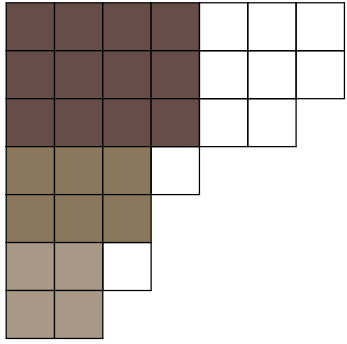
 \label{fig-1}
$$ \begin{ytableau}
   *(Beaver) & *(Beaver) &*(Beaver)& *(Beaver)& & &   \\
  *(Beaver) & *(Beaver) & *(Beaver)&*(Beaver) & & & \\
     *(Beaver)  & *(Beaver) & *(Beaver)& *(Beaver)& & \\
      *(Shadow)& *(Shadow) &*(Shadow) & \\
       *(Shadow)  &   *( Shadow)  & *( Shadow)     \\
       *(grullo) &     *(grullo)  &  \\
       *(grullo)   &   *(grullo)
  \end{ytableau}$$
\caption{ The Young diagram of the partition $\lambda=(7,7,6,4,3,3,2)$. 
 Durfee rectangles and Durfee squares in the figure show that  $\lambda \in \mathcal{D}_{4,2}(32)$.}
\end{figure}

The set $\mathcal{D}_{r,i}(n)$ contains the partitions of $n$ with \textit{at most} $r-i$ successive horizontal Durfee rectangles followed by $i-1$ Durfee squares and no parts after the last square/rectangle (when we have the parts equal to one we can consider the horizontal Durfee rectangles of size $1 \times 0$).  For instance the partition $\lambda=(7,7,6,4,3,3,2)$ belongs to $\mathcal{D}_{4,2}(32)$  (see 
Figure \ref{fig-1}).  Notice also that  $ \lambda \in \mathcal{D}_{6,4}(32)$. 
\bigskip

 By equation ($2.4$) of \cite{ADJM1} the generating series of $\mathcal{D}_{r,i}(n)$ is equal to:
\begin{equation}\label{D}
\sum_{n\in \mathbb{N}} \mathcal{D}_{r,i}(n) q^n= \sum_{\partial_1\geq \cdots \geq \partial_{r-1}\geq 0} \frac{q^{\partial_1^2+\cdots+\partial_{r-1}^2-\partial_{1}-\cdots-\partial_{r-i}}}{(q)_{\partial_1-\partial_2}\cdots (q)_{\partial_{r-2}-\partial_{r-1}}(q)_{\partial_{r-1}}}(1-q^{\partial_{r-i}}),
\end{equation}
where $\partial_j$ is the \textit{largest} side of the $j$-th Durfee square/rectangle and for $i=r$ we have:
 $$\sum_{n\in \mathbb{N}} \mathcal{D}_{r,r}(n) q^n= \sum_{\partial_1\geq \cdots \geq \partial_{r-1}\geq 0} \frac{q^{\partial_1^2+\cdots+\partial_{r-1}^2}}{(q)_{\partial_1-\partial_2}\cdots (q)_{\partial_{r-2}-\partial_{r-1}}(q)_{\partial_{r-1}}}.$$
 
 For instance in the previous example we have $\partial_1=4, \partial_2=3$ and $\partial_3=2$  (see 
Figure \ref{fig-1}).
 \bigskip
 
  Denote by $b_{r,i}(m,n)$ (respectively by $c_{r,i}(m,n)$ and by $d_{r,i}(m,n)$) the number of partitions of $\mathcal{B}_{r,i}(n)$ (respectively of $\mathcal{C}_{r,i}(n)$ and of $\mathcal{D}_{r,i}(n)$) with exactly $m$ parts. Based on the results of \cite{BMS} and \cite{PE}, we can deduce the following:
  
  \begin{theo} \label{refinement}Let $m,n$ be integers. We have the identity
  \[b_{r,i}(m,n)=c_{r,i}(m,n)\]
   \end{theo}
  
  \begin{proof}
We prove the theorem in the case $i=r;$ the other cases follow in a similar way. 
Consider the differential ideal  $\mathcal{I}$ of the ring
 \[
  \S_1 := \mathbf{K}[x_i,~i\in \mathbf{Z}_{\geq 1}],\]
generated by $x_1^r$ together with all its iterated derivatives $D^j(x_1^r)$ for $j\in \mathbf{Z}_{\geq 1}$. Here, as in the introduction, $D$ denotes the derivation defined by $D(x_i)=x_{i+1}$ and extended to the whole algebra via the Leibniz rule. The iterated derivatives are defined recursively by $D^j(x_1^r)=D(D^{j-1}(x_1^r))$ with $D^1=D$ and $D^0$ the identity. As explained in the introduction, the ideal $\mathcal{I}$ is weighted homogeneous if we assign to $x_i$ the weight $i,$ but it is also homogeneous with respect to the usual degree, \emph{i.e.}, if we assign to $x_i$ the weight $1$ for  every $i$ (this is follows from the fact that $x_1^r$ is homogeneous and hence all its iterative derivatives are). We deduce that ring 
\[
\mathcal{R}:=  \S_1 / \mathcal{I}   
=\sum_{i\in \mathbf{Z}_{\geq 0}, d\in \mathbf{Z}_{\geq 0}}\mathcal{R}_{i,d}\]
is bi-graded, where $\mathcal{R}_{i,d}$ is the bi-homogeneous component of elements of weight $i$ and degree $d$.
Consider the partial monomial ordering on ${\S}_1$
 defined  by the usual degree, and then refine it either by the weighted reverse lexicographical ordering which gives the total monomial ordering $<_1$ or by the weighted lexicographical ordering $<_2$ (here, the weight of $x_i$ is 
$i$). That is,  given  two monomials $M,M'\in  \S_1$  then 
$M<_1 M'$ (respectively $M<_2 M'$) if the degree of $M$ is strictly smaller than the degree of $M'$ or they are of equal degrees and $M$ is smaller than $M'$ with respect to the weighted reverse lexicographical ordering (resp. the weighted lexicographical ordering). Now, if one begins with the homogeneous generators of $\mathcal{I}$ defined above, and compute the $S-$polynomial (with respect to $<_1$ -respectively $<_2$-, as in the Buchberger algorithm) of two generators, we see  that this is the same thing as computing the $S-$polynomial of the same two generators with respect to the weighted reverse lexicographical ordering (resp. the
the weighted lexicographical ordering); \textit{i.e}, the fact that we consider the partial order with respect to the usual degree does not change the algorithm. We deduce that the initial ideal with respect to the weighted reverse lexicographical (respectively the weighted lexicographical ordering) is equal to the initial ideal with respect to the ordering $<_1$ (respectively $<_2).$ Let us call
 $\mathcal{I}'$
 (respectively 
$\mathcal{I}''$)
 the initial ideal of $\mathcal{I}$ with respect to the weighted reverse lexicographical ordering (resp. with respect to the weighted lexicographical ordering). The rings 
\[\mathcal{R}':=  \S_1 / \mathcal{I}'
=\sum_{i\in \mathbf{Z}_{\geq 0}, d\in \mathbf{Z}_{\geq 0}}\mathcal{R}'_{i,d} \] 
\[\mathcal{R}'':= \S_1/ \mathcal{I}''
=\sum_{i\in \mathbf{Z}_{\geq 0}, d\in \mathbf{Z}_{\geq 0}}\mathcal{R}''_{i,d} \]
are again bi-graded where $i$ indicates the weight and $d$ the degree. Thanks to the fact that the ideals 
 $\mathcal{I}'$ and $\mathcal{I}''$
are initial ideals with respect to  monomial orderings that respect the weight and the degree, we have  the equalities:
\[\mbox{dim}_{\mathbf{K}}~\mathcal{R}'_{i,d}=\mbox{dim}_{\mathbf{K}}~\mathcal{R}_{i,d}~~~~~~\mbox{and}~~~~~~\mbox{dim}_{\mathbf{K}}~\mathcal{R}''_{i,d}=\mbox{dim}_{\mathbf{K}}~\mathcal{R}_{i,d}.\] 
It follows from \cite{BMS} that 
\[\mbox{dim}_{\mathbf{K}}~\mathcal{R}'_{i,d}=b_{r,i}(m,n)
   \]
and from \cite{Af} and \cite{PE} that 
\[\mbox{dim}_{\mathbf{K}}~\mathcal{R}''_{i,d}=c_{r,i}(m,n).\] Hence the equality $b_{r,i}(m,n)=c_{r,i}(m,n)$ holds.
  \end{proof}

 We have the following proposition in which we give the generating series of the partitions verifying the difference condition of Gordon's identities \emph{with a fixed length $m\geq0$}:

 \begin{prop}\label{m} For the integers $m,n \geq 0, r \geq 2$ and $1\leq i \leq r$  we have:
 $$\sum_{n\in \mathbb{N}} b_{r,i}(m,n) q^n= \sum_{\underset{\partial_1\geq \cdots \geq \partial_{r-1}\geq 0}{\partial_1+\cdots+\partial_{r-1}=m+r-i}} \frac{q^{\partial_1^2+\cdots+\partial_{r-1}^2-\partial_{1}-\cdots-\partial_{r-i}}(1-q^{\partial_{r-i}})}{(q)_{\partial_1-\partial_2}\cdots (q)_{\partial_{r-2}-\partial_{r-1}}(q)_{\partial_{r-1}}}$$
 \end{prop} 
 
 \begin{proof} 
By Theorem~\ref{refinement}, we have $b_{r,i}(m,n)=c_{r,i}(m,n)$. By Theorem $2.3$ in \cite{ADJM1}, we have $c_{r,i}(m,n)=d_{r,i}(m,n),$ hence the three numbers are equal. If we denote by $\partial_j$ he largest side of the $j$-th Durfee square/rectangle, then we can compute the generating series of $d_{r,i}(m,n)$  in the same way as in the computations in \cite{ADJM1} of the generating series of $D_{r,i}(n)$. The difference is that here, the length of the partitions must be set to $m,$ which is equal to $(\partial_1-1)+\cdots+(\partial_{r-i}-1)+\partial_{r-i+1}+\cdots+\partial_{r-1}.$ We obtain the same formula as in Equation (\ref{D}) with the fixed length $m:$
$$
\sum_{n\in \mathbb{N}} d_{r,i}(m,n) q^n= \sum_{\underset{\partial_1\geq \cdots \geq \partial_{r-1}\geq 0}{(\partial_1-1)+\cdots+(\partial_{r-i}-1)+\partial_{r-i+1}+\cdots+\partial_{r-1}=m}} \frac{q^{\partial_1^2+\cdots+\partial_{r-1}^2-\partial_{1}-\cdots-\partial_{r-i}}(1-q^{\partial_{r-i}})}{(q)_{\partial_1-\partial_2}\cdots (q)_{\partial_{r-2}-\partial_{r-1}}(q)_{\partial_{r-1}}}.
$$
\end{proof} 
We are now ready to determine the generating series of $G_{r,\ell}(n):$
\begin{prop}\label{G} Let $\ell, n \geq 0$ and $r\geq 2$ be the integers. The generating series $G_{r,\ell}$ of $G_{r,\ell}(n)$ is equal to:

$$\sum_{\underset{\partial_1\geq \cdots \geq\partial_{r-1}\geq 0}{\partial_1+\cdots+\partial_{r-1}=0}}^{\ell+r-1} \frac{q^{\partial_1^2+\cdots+\partial_{r-1}^2}}{(q)_{\partial_1-\partial_2}\cdots (q)_{\partial_{r-2}-\partial_{r-1}}(q)_{\partial_{r-1}}}-\sum_{\underset{\partial_1\geq \cdots \geq \partial_{r-1}\geq 0}{\partial_1+\cdots+\partial_{r-1}=\ell+r}} \frac{q^{\partial_1^2+\cdots+\partial_{r-1}^2-(\ell+r)}(1-q^{\ell+r})}{(q)_{\partial_1-\partial_2}\cdots (q)_{\partial_{r-2}-\partial_{r-1}}(q)_{\partial_{r-1}}}.$$
\end{prop}
\begin{proof}Let $\lambda\in \mathcal{G}_{r,\ell}(n).$ Therefore its parts verify  Gordon's difference condition and we have $$L(\lambda)-\#\{\text{parts equal to $1$ of $\lambda$}\}\leq \ell.$$

Since $\lambda$ has at most $r-1$ parts equal to one, so $\mathcal{G}_{r,\ell}(n)$ is equal to the union of the sets of partitions $\lambda$ of $\mathcal{B}_{r,r}(n)$ with exactly $j$ parts equal to one and with $L(\lambda)-j\leq \ell,$ for $j= 0,\ldots,r-1.$
\\


Recall that $b_{r,i}(m,n)$ is the number of partitions of $n$ satisfying the Gordon's difference condition with exactly $m$ parts and with at most $i-1$ parts equal to $1$. Set $b_{r,0}(m,n)=0$. Note that
$b_{r,i}(m,n)-b_{r,i-1}(m,n)$
is equal to the number of partitions of $\mathcal{B}_{r,r}(n)$  with exactly $m$ parts among which exactly $i-1$ parts are equal to one (see the proof of Theorem  $7.5$ from \cite{A}). Thus 
$$\mathcal{G}_{r,\ell}(n)=\bigcup_{i=1}^{r} \{\lambda \in \mathcal{B}_{r,r}(n)| \ m-(i-1)\leq \ell \ \text{ and $\lambda $ is counted by} \ b_{r,i}(m,n)-b_{r,i-1}(m,n)\},$$
and since the sets in the union are disjoint, we have:
$$
G_{r,\ell}(n) = \sum_{i=1}^{r} \sum_{m=0}^{\ell+i-1} \Big( b_{r,i}(m,n)-b_{r,i-1}(m,n)\Big)=\sum_{m=0}^{\ell+r-1} b_{r,r}(m,n)-\sum_{i=1}^{r-1} b_{r,i}(\ell+i,n).
$$

Considering the generating series of $G_{r,\ell}(n)$ using the above equation and then applying Proposition \ref{m} on $b_{r,r}(m,n)$ and on $b_{r,i}(\ell+i)$ yields:

$$G_{r,\ell}=\sum_{n\in \mathbb{N}} G_{r,\ell}(n) q^n =   \sum_{\underset{\partial_1\geq \cdots \geq\partial_{r-1}\geq 0}{\partial_1+\cdots+\partial_{r-1}=0}}^{\ell+r-1} \frac{q^{\partial_1^2+\cdots+\partial_{r-1}^2}}{(q)_{\partial_1-\partial_2}\cdots (q)_{\partial_{r-2}-\partial_{r-1}}(q)_{\partial_{r-1}}}$$
 
 $$ -\sum_{i=1}^{r-1} \sum_{\underset{\partial_1\geq \cdots \geq \partial_{r-1}\geq 0}{(\partial_1-1)+\cdots+(\partial_{r-i}-1)+\partial_{r-i+1}+\cdots+\partial_{r-1}=\ell+i}} \frac{q^{\partial_1^2+\cdots+\partial_{r-1}^2-\partial_1-\cdots-\partial_{r-i}}}{(q)_{\partial_1-\partial_2}\cdots (q)_{\partial_{r-2}-\partial_{r-1}}(q)_{\partial_{r-1}}}(1-q^{\partial_{r-i}}).$$

Note that for each $1\leq i \leq r$ the condition of the last sum is that $\partial_1+\cdots+\partial_{r-1}=\ell+r.$ Therefore we can change the order of the last two sums. Developing now the last one and then putting the fractions all over the same denominator, we obtain:

\begin{align*}
\sum_{n\in \mathbb{N}} G_{r,\ell}(n) q^n &=   \sum_{\underset{\partial_1\geq \cdots \geq\partial_{r-1}\geq 0}{\partial_1+\cdots+\partial_{r-1}=0}}^{\ell+r-1} \frac{q^{\partial_1^2+\cdots+\partial_{r-1}^2}}{(q)_{\partial_1-\partial_2}\cdots (q)_{\partial_{r-2}-\partial_{r-1}}(q)_{\partial_{r-1}}}
 \\
  &- \sum_{\underset{\partial_1\geq \cdots \geq \partial_{r-1}\geq 0}{\partial_1+\cdots+\partial_{r-1}=\ell+r}} \frac{q^{\partial_1^2+\cdots+\partial_{r-1}^2}\Big(q^{-\partial_1-\cdots-\partial_{r-1}}(1-q^{\partial_{r-1}})+\cdots+q^{-\partial_1}(1-q^{\partial_{1}})\Big)}{(q)_{\partial_1-\partial_2}\cdots (q)_{\partial_{r-2}-\partial_{r-1}}(q)_{\partial_{r-1}}}\\  
  &=   \sum_{\underset{\partial_1\geq \cdots \geq\partial_{r-1}\geq 0}{\partial_1+\cdots+\partial_{r-1}=0}}^{\ell+r-1} \frac{q^{\partial_1^2+\cdots+\partial_{r-1}^2}}{(q)_{\partial_1-\partial_2}\cdots (q)_{\partial_{r-2}-\partial_{r-1}}(q)_{\partial_{r-1}}}
 \\
  &- \sum_{\underset{\partial_1\geq \cdots \geq \partial_{r-1}\geq 0}{\partial_1+\cdots+\partial_{r-1}=\ell+r}} \frac{q^{\partial_1^2+\cdots+\partial_{r-1}^2}\Big(q^{-\partial_1-\cdots-\partial_{r-1}}-1\Big)}{(q)_{\partial_1-\partial_2}\cdots (q)_{\partial_{r-2}-\partial_{r-1}}(q)_{\partial_{r-1}}}\\
\end{align*}

Replacing $-\partial_1-\cdots-\partial_{r-1}$ by $-(\ell+r)$ and factorizing by $q^{-(\ell+r)}$ in the 
numerator
of the last sum give the result.

\end{proof}

\begin{rem}\label{lamdag} For each partition $\lambda \in \mathcal{F}_{r}(n)$ of type $4.b$ we have $\lambda_g\in \mathcal{G}_{r,k+i_k-r}(n)$ for some positive integer $n,$ where $k$ is the smallest part of $\lambda_b$ and $i_k$ is the $k$-th smallest part of $\lambda_r.$ 
\end{rem}

\begin{rem}\label{Hil} Note that the partitions of $\mathcal{G}_{r,\ell}(n)$ are those associated with the monomials of the algebra $$\frac{\mathbf{K}[z_1,z_2,\cdots]}{(I_r,z_{j_1}\cdots z_{j_{\ell+1}}| \ 2\leq j_1 \leq \cdots \leq j_{\ell+1})}
.$$
 Thus, the Hilbert-Poincar\'e series of this algebra is equal to the generating series of $G_{r,\ell}(n)$ which is represented in the Proposition \ref{G}.
\end{rem}

For a positive integer $m\geq 0$ denote by 
\[ X_m:=\sum_{n\in \mathbb{N}}d_{r,r}(m,n) q^n .\] 
By Proposition \ref{m} we have the equality
$$X_m = \sum_{\underset{\partial_1\geq \cdots \geq\partial_{r-1}\geq 0}{\partial_1+\cdots+\partial_{r-1}=m}} \frac{q^{\partial_1^2+\cdots+\partial_{r-1}^2}}{(q)_{\partial_1-\partial_2}\cdots (q)_{\partial_{r-2}-\partial_{r-1}}(q)_{\partial_{r-1}}}.$$ 
Thus with this notation the first sum in $G_{r,\ell}$ is equal to $\sum_{m=0}^{\ell+r-1}X_m.$
\\

We have the following lemma about this part of $G_{r,\ell}$ which we will use later in the last section in order to compute the generating series of the partitions of $\mathcal{F}_r$ of type $4.b$ (see Lemma \ref{S4ii} below):
\begin{lem}\label{X} The following equation holds:
$$\sum_{\ell\geq 0} q^{r+\ell}(q)_{r+\ell-1} \sum_{m=0}^{\ell+r-1} X_m = 1-\frac{\mathbf{G}_r}{\mathbf{H}}+\sum_{m \geq r }X_m (q)_m.$$
\end{lem}
\begin{proof}
We have:
\begin{align*}
\sum_{\ell\geq 0} q^{r+\ell}(q)_{r+\ell-1} \sum_{m=0}^{\ell+r-1} X_m &=  q^r (q)_{r-1}  \sum_{m=0}^{r-1} X_m 
 \\
  &+ q^{r+1} (q)_{r} ( \sum_{m=0}^{r-1} X_m+X_r) \\
  &+ q^{r+2} (q)_{r+1} ( \sum_{m=0}^{r-1} X_m+X_r+ X_{r+1})
 \\
 & \vdots
\end{align*}
Summing over the columns we obtain:

 $$\sum_{m=0}^{r-1} X_m \sum_{k\geq r} q^{k} (q)_{k-1} 
 + X_r \sum_{k\geq r+1} q^{k} (q)_{k-1} 
  +  X_{r+1} \sum_{k\geq r+2} q^{k} (q)_{k-1} 
 +\cdots$$
 $$=\sum_{m=0}^{r-1} X_m \sum_{k\geq r} q^{k} (q)_{k-1}+ \sum_{m\geq r} X_m \sum_{k\geq m+1} q^{k} (q)_{k-1}.$$
 Using twice Lemma \ref{H} yields:
 $$\sum_{m=0}^{r-1} X_m \Big((q)_{r-1}-\frac{1}{\mathbf{H}}\Big)+ \sum_{m\geq r} X_m \Big((q)_{m}-\frac{1}{\mathbf{H}}\Big)= (q)_{r-1}\sum_{m=0}^{r-1} X_m- \frac{1}{\mathbf{H}} \sum_{m\geq 0} X_m+\sum_{m\geq r} X_m (q)_{m}.$$
 Recall that:
 $$\sum_{m\geq 0} X_m=\sum_{n\in \mathbb{N}} \sum_{m\geq0} d_{r,r}(m,n)q^n=\sum_{n\in \mathbb{N}} D_{r,r}(n)q^n=\mathbf{G}_r,$$
 and that:
 $$\sum_{m=0}^{r-1} X_m=\sum_{n\in \mathbb{N}} \sum_{m=0}^{r-1} d_{r,r}(m,n)q^n.$$
 This is the generating function of the partitions in $\mathcal{D}_{r,r}(n)$ with \textit{at most $r-1$ parts}. Note that each partition of length $\leq r-1$ has also \textit{at most $r-1$ successive Durfee squares} and belongs to $\mathcal{D}_{r,r}(n)$ for some positive integer $n.$ This means that $ \sum_{m=0}^{r-1} d_{r,r}(m,n)$ is equal to the number of partitions of $n$ with at most $r-1$ parts. Thus by Equation (\ref{red1}) we have:
 $$\sum_{m=0}^{r-1} X_m=\underbrace{1}_{\text{empty set}}+P_{\leq r-1}=1+\frac{1}{(q)_{r-1}}-1=\frac{1}{(q)_{r-1}}.$$
 So we have:
 $$\sum_{\ell\geq 0} q^{r+\ell}(q)_{r+\ell-1} \sum_{m=0}^{\ell+r-1} X_m = 1-\frac{\mathbf{G}_r}{\mathbf{H}}+\sum_{m \geq r }X_m (q)_m.$$
\end{proof}

\section{Hilbert-Poincar\'e series of the algebra $ \S/J_r$}
As we mentioned in the second section, the Hilbert-Poincar\'e series of $ \S/J_r$ is equal to the generating series of $\mathcal{F}_{r}(n)$. We have also computed the generating series of partitions of type $1$ and type $2$ of $\mathcal{F}_{r}(n)$. In this section we compute the generating series of its other two types and finally we give the generating series of $\mathcal{F}_{r}(n).$
 \begin{lem}\label{S3} For the integers $r\geq 2$ and $n\geq 0$ the generating series  of the partitions of type $3$ of $\mathcal{F}_{r}(n)$ is equal to:
 $$S_3:=1-2\mathbf{H}+\mathbf{H}^2\ \frac{(q)_{r-1}}{(q)_1}.$$
 \end{lem}
 \begin{proof}
 In order to give $S_3$
we will  compute the generating series of the partitions in $\mathcal{F}_{r}(n)$ of each type $3.a$ and $3.b$ and then add them together.
 Note that each partition $\lambda \in \mathcal{F}_{r}(n)$ of type $3$ has \textit{exactly} two colors: black and red. 

 For a partition $\lambda$ of type $3.a,$ denote the smallest part of $\lambda_b$ by $k$. Recall that in this case the length $\ell_r$ of $\lambda_r$ is between $1$ and $k-1$ (thus $k\geq 2$). Note that the partitions of type $3.a$ are generated by the same $q$-series as the generating series of the $1$-colored partitions $\lambda_b$ multiplied by the generating series of the $1$-colored partitions $\lambda_r.$ Thus, by using Lemma \ref{black} and equation (\ref{red1}), the partitions of type $3.a$ are generated by:
$$
 \sum_{k\geq 2} \underbrace{\frac{q^k}{\prod_{j\geq k}(1-q^j)}}_{\text{Black sub-partitions of type $3.a$}}\Big(\underbrace{\frac{1}{(q)_{k-1}}-1}_{\text{red sub-partitions of type $3.a$}}\Big)=\mathbf{H}\frac{q^2}{(q)_1}-\sum_{k\geq2} \frac{q^k}{\prod_{j\geq k}(1-q^j)}.
$$
Note that by Lemma \ref{black} the sum $\sum_{k\geq2} \frac{q^k}{\prod_{j\geq k}(1-q^j)}$ generates all partitions except the empty set and the partitions whose smallest part is equal to $1$. Therefore, this sum is equal to $\mathbf{H}-1-q \mathbf{H}$ and so the generating series of the partitions $3.a$ is equal to:
\begin{equation}\label{3i}
S_{3,a}:=\mathbf{H} \ \Big(\frac{2q-1}{(q)_1}\Big)+1.
\end{equation}

 For the partitions $\lambda$ of type $3.b,$ denote the smallest part of $\lambda_b$ by $k.$ Recall that the length $\ell_r$ of $\lambda_r$ is greater than or equal to $k$ with $k$-th smallest part $i_k\geq 1$ such that $k+i_k\geq r.$ Thus, by using Lemma \ref{black} and Lemma \ref{red2}, we obtain the generating series of the partitions of type $3.b$ as follows:
 
\begin{equation}\label{3.b}
S_{3,b}=\sum_{k\geq 1} \sum_{\underset{i_k\geq 1}{i_k\geq r-k}}\underbrace{\frac{q^k}{\prod_{j\geq k}(1-q^j)}}_{\text{Black sub-partitions of type $3.b$}}\Big(\underbrace{\mathbf{H}. q^{k+i_k-1} \frac{(q)_{k+i_k-2}}{(q)_{k-1}}}_{\text{red sub-partitions of type $3.b$}}\Big) .
\end{equation}

This is equal to:

\begin{align*}
\mathbf{H}^2 \sum_{k\geq 1}  q^k \sum_{\underset{i_k\geq 1}{i_k\geq r-k}} q^{k+i_k-1} (q)_{k+i_k-2}=   \mathbf{H}^2 & \Big[ \sum_{k=1}^{r-1}  q^k \sum_{{i_k\geq r-k}} q^{k+i_k-1} (q)_{k+i_k-2}  
 \\
  & +\sum_{k\geq r} q^k \sum_{i_k\geq 1} q^{k+i_k-1} (q)_{k+i_k-2}   \Big],\\
\end{align*}
applying twice Lemma \ref{H}, once for $m=r-1$ and once for $m=k$, we obtain:
\begin{align*}
   \mathbf{H}^2 & \Big[ \sum_{k=1}^{r-1}  q^k \big((q)_{r-2}-\frac{1}{\mathbf{H}}\big) +\sum_{k\geq r} q^k  \big((q)_{k-1}-\frac{1}{\mathbf{H}}\big)   \Big]
 \\
  = \mathbf{H}^2 & \Big[  \big((q)_{r-2}-\frac{1}{\mathbf{H}}\big) \frac{q(1-q^{r-1})}{(q)_1} +\sum_{k\geq r} q^k (q)_{k-1}-\frac{1}{\mathbf{H}}\frac{q^r}{(q)_1}   \Big].\\
\end{align*}
Applying once again Lemma \ref{H} for $m=r$ and simplifying we get the generating series of the partitions of type $3.b$ as follows:
$$S_{3,b}=\mathbf{H}^2 \frac{(q)_{r-1}}{(q)_1}-\mathbf{H}\frac{1}{(q)_1},$$

Adding this last formula to the formula (\ref{3i}), yields the result.
 \end{proof}
 
\begin{lem}\label{S4i} For the integers $r\geq 2$ and $n\geq 0$ the generating series of the partitions of type $4.a$ of $\mathcal{E}_{r-1}(n)$ is equal to:
 $$S_{4,a}:=(\mathbf{G}_{r}-1)S_{3,a}.$$
 \end{lem} 
 \begin{proof}
 Note that for each partition $\lambda$ of type $4.a$ the sub-partitions $\lambda_b, \lambda_r$ and $\lambda_g$ are all non-empty. Denote the smallest part of $\lambda_b$ by $k$. Recall that in this case 
  $\ell_r\leq k-1$ and that $\lambda_g \in \mathcal{B}_{r,r}(n).$ Thus, $\lambda_g$ is generated by $\mathbf{G}_r-1$ and using Lemma \ref{black} together with equation (\ref{red1}) we have:
 
$$S_{4,a}= \sum_{k\geq 2} \underbrace{\frac{q^k}{\prod_{j\geq k}(1-q^j)}}_{\text{Black sub-partitions of type $4.a$}}\Big(\underbrace{\frac{1}{(q)_{k}}-1}_{\text{red sub-partitions of type $4.a$}}\Big) \Big(\underbrace{\mathbf{G}_r-1}_{\text{Green sub-partitions of type $4.a$}}\Big).$$
We can take $\mathbf{G}_r-1$ out of the sum. The remaining part is equal to $S_{3,a}$ (see the proof of Lemma \ref{S3}). 
 \end{proof}
 
 \begin{prop}\label{S4ii} For the integers $r\geq 2$ and $n\geq 0$ the generating series of the partitions of type $4.b$ of $\mathcal{F}_{r}(n)$ is equal to:
 $$S_{4,b}:=\frac{\mathbf{H}^2-\mathbf{H} \, \mathbf{G}_r}{(q)_1}-S_{3,b}.$$
 \end{prop} 
 
 \begin{proof} For each partition $\lambda$ of $\mathcal{F}_{r}$ of type $4.b$ denote by $k$ the smallest part of $\lambda_b.$ Recall that in this case $\ell_r\geq k$  with $k+i_k \geq r$ and $\ell_g-\#\{1_g\}<k+i_k-r+1,$ where $i_k$ is the $k$-th smallest part  of $\lambda_r.$ Recall also that $\lambda_g \in \mathcal{G}_{r,k+i_k-r}(n)$ for some positive integer $n$ (see Remark \ref{lamdag}). Thus using Lemma \ref{black} and Lemma \ref{red2}  we have::
 
 $$ \sum_{k\geq 1}\sum_{\underset{i_k\geq 1}{i_k\geq r-k}} \underbrace{\frac{q^k}{\prod_{j\geq k}(1-q^j)}}_{\text{Black sub-partitions of type $4.b$}}\Big(\underbrace{\mathbf{H}. q^{k+i_k-1} \frac{(q)_{k+i_k-2}}{(q)_{k-1}}}_{\text{red sub-partitions of type $4.b$}}\Big) \Big(\underbrace{
 G_{r,k+i_k-r}-1}_{\text{Green sub-partitions of type $4.b$}}\Big).$$
 By simplifications and using equation (\ref{3.b}) we obtain:
 $$S_{4.b}=\mathbf{H}^2 \sum_{k\geq 1}\sum_{\underset{i_k\geq 1}{i_k\geq r-k}} q^{2k+i_k-1} (q)_{k+i_k-2} G_{r,k+i_k-r}-S_{3,b}.$$
 Thus we have:
 \begin{align*}
S_{4.b}+S_{3.b}=   \mathbf{H}^2 & \Big[ \sum_{k=1}^{r-1} \sum_{{i_k\geq r-k}} q^{2k+i_k-1} (q)_{k+i_k-2} G_{r,k+i_k-r}
 \\
  & +\sum_{k\geq r} \sum_{i_k\geq 1} q^{2k+i_k-1} (q)_{k+i_k-2} G_{r,k+i_k-r}  \Big]\\
\end{align*}
 Developing the second sum of each double sum gives:
 \begin{align*}
 \mathbf{H}^2 & \Big[ \sum_{k=1}^{r-1}\Big( q^{r+k-1} (q)_{r-2} G_{r,0}+ q^{r+k} (q)_{r-1} G_{r,1}+ q^{r+k+1} (q)_{r} G_{r,2}+\cdots \Big)
 \\
  & +\sum_{k\geq r} \Big( q^{2k} (q)_{k-1} G_{r,k-r+1}+q^{2k+1} (q)_{k} G_{r,k-r+2}+ q^{2k+2} (q)_{k+1} G_{r,k-r+3}+\cdots \Big) \Big].\\
 \end{align*}
 Developing now these two sums we obtain $\mathbf{H}^2 $ times:

  \xymatrix{
      q^{r} (q)_{r-2} G_{r,0} & +q^{r+1} (q)_{r-1}G_{r,1}&+q^{r+2} (q)_{r} G_{r,2}+\cdots \\
       +q^{r+1} (q)_{r-2} G_{r,0}&+q^{r+2} (q)_{r-1}G_{r,1}&+q^{r+3} (q)_{r} G_{r,2}+\cdots  \\
   \vdots   &  \vdots & \vdots \\
     +q^{2r-2} (q)_{r-2} G_{r,0}&+q^{2r-1} (q)_{r-1}G_{r,1}&+q^{2r} (q)_{r} G_{r,2}+\cdots \\     
     &+ q^{2r} (q)_{r-1} G_{r,1}&+q^{2r+1} (q)_{r}G_{r,2}+\cdots\\
     & &+  q^{2r+2} (q)_{r} G_{r,2}+\cdots\\
     & & \ddots
 }
  Summing over the columns we get:
   \begin{align*}
 S_{4.b}+S_{3.b} & =   \mathbf{H}^2 \sum_{\ell\geq 0} (q^{r+\ell}+\cdots+q^{2r+2\ell-2})(q)_{r+\ell-2}G_{r,\ell}
 \\
  & =\mathbf{H}^2 \sum_{\ell\geq 0} \frac{q^{r+\ell}(1-q^{r+\ell-1})}{(q)_1}(q)_{r+\ell-2}G_{r,\ell}\\
  \\
  & =\mathbf{H}^2 \sum_{\ell\geq 0} \frac{q^{r+\ell}}{(q)_1}(q)_{r+\ell-1}G_{r,\ell}.
 \end{align*}
  Using Proposition \ref{G}, we now replace $G_{r,\ell}$:
  \begin{align*}
 S_{4.b}+S_{3.b} =\frac{\mathbf{H}^2}{(q)_1} \sum_{\ell\geq 0} q^{r+\ell}(q)_{r+\ell-1} & \Big( \sum_{\underset{\partial_1\geq \cdots \geq\partial_{r-1}\geq 0}{\partial_1+\cdots+\partial_{r-1}=0}}^{\ell+r-1} \frac{q^{\partial_1^2+\cdots+\partial_{r-1}^2}}{(q)_{\partial_1-\partial_2}\cdots (q)_{\partial_{r-2}-\partial_{r-1}}(q)_{\partial_{r-1}}}\\
 &-\sum_{\underset{\partial_1\geq \cdots \geq \partial_{r-1}\geq 0}{\partial_1+\cdots+\partial_{r-1}=\ell+r}} \frac{q^{\partial_1^2+\cdots+\partial_{r-1}^2-(\ell+r)}(1-q^{\ell+r})}{(q)_{\partial_1-\partial_2}\cdots (q)_{\partial_{r-2}-\partial_{r-1}}(q)_{\partial_{r-1}}}\Big). \\
 \end{align*}
 Simplifying and using the notations of Lemma \ref{X}, we have:
  \begin{align*}S_{4.b}+S_{3.b}=\frac{\mathbf{H}^2}{(q)_1} & \Big( \sum_{\ell\geq 0} q^{r+\ell}(q)_{r+\ell-1} \sum_{m=0}^{\ell+r-1} X_m\\&-\sum_{\ell\geq 0} (q)_{r+\ell} \sum_{\underset{\partial_1\geq \cdots \geq \partial_{r-1}\geq 0}{\partial_1+\cdots+\partial_{r-1}=\ell+r}} \frac{q^{\partial_1^2+\cdots+\partial_{r-1}^2}}{(q)_{\partial_1-\partial_2}\cdots (q)_{\partial_{r-2}-\partial_{r-1}}(q)_{\partial_{r-1}}}\Big).
   \end{align*}
 Using Lemma \ref{X} and changing the variable $\ell+r$ by $m$ in the last two sums, give:
 
 \begin{align*}S_{4.b}+S_{3.b}=\frac{\mathbf{H}^2}{(q)_1}& \Big(1-\frac{\mathbf{G}_r}{\mathbf{H}}+\sum_{m \geq r }X_m (q)_m\\-& \sum_{m\geq r}  (q)_{m} \sum_{\underset{\partial_1\geq \cdots \geq \partial_{r-1}\geq 0}{\partial_1+\cdots+\partial_{r-1}=m}} \frac{q^{\partial_1^2+\cdots+\partial_{r-1}^2}}{(q)_{\partial_1-\partial_2}\cdots (q)_{\partial_{r-2}-\partial_{r-1}}(q)_{\partial_{r-1}}}\Big). 
 \end{align*}
which is equal to:
$$\frac{\mathbf{H}^2}{(q)_1} \Big(1-\frac{\mathbf{G}_r}{\mathbf{H}}+\sum_{m \geq r }X_m (q)_m -\sum_{m \geq r }X_m (q)_m \Big)=\frac{\mathbf{H}^2-\mathbf{H} \, \mathbf{G}_r}{(q)_1}. 
$$
 \end{proof}
 Recall that $F_r(n)$ is the number of partitions of weight $n$ associated to the algebra $\S/J_r$. We are now ready to  prove 
  the main result of this paper, stated as Theorem \ref{th-frn} at the introduction,
  and which  gives us a family of an infinite number of $3$-colored partition identities associated with $A_{r-1}$:
 \begin{theo}\label{main} For all $n \in \n$ and all integers $r \geq 2$, the numbers $F_r(n)$ are equal:
\[
F_2(n) = F_3(n) = F_4(n) = \cdots.
\]
Moreover, this common value equals the number of integer partitions of $n$ in which the part $1$ may appear in $3$ colors, while all other positive integers may appear in $2$ colors.
 
 \end{theo}
 
 \begin{proof}
As we have mentioned before, for any integer $r\geq 2$, the generating series of $F_r(n)$ is equal to the sum of the generating series of the partitions of type $1$ to $4$ of $\mathcal{F}_r(n).$ By Propositions \ref{S4i} and \ref{S4ii}, we have:

\begin{align*}\sum_{n\in \mathbb{N}} F_r(n)q^n&=S_1+S_2+S_3+S_{4,a}+S_{4,b}\\
 & =S_1+S_2+S_3+(\mathbf{G}_{r}-1)S_{3,a}+\frac{\mathbf{H}^2-\mathbf{H} \, \mathbf{G}_r}{(q)_1}-S_{3,b}\\
 & =S_1+S_2+\mathbf{G}_{r} S_{3,a}+\frac{\mathbf{H}^2-\mathbf{H}\mathbf{G}_r}{(q)_1}.\\
  \end{align*}
   Using Equations (\ref{S1}) and (\ref{S2}) together with Equation (\ref{3i}) gives:
 \begin{align*}\sum_{n\in \mathbb{N}} F_r(n)q^n
  &=\Big(2 \mathbf{H}+\mathbf{G}_r-2\Big)+\Big(2(\mathbf{H}-1)(\mathbf{G}_r-1)\Big)\\
  &+\mathbf{G}_{r}\Big(\mathbf{H} \ \Big(\frac{2q-1}{(q)_1}\Big)+1\Big)+\frac{\mathbf{H}^2-\mathbf{H} \,\mathbf{G}_r}{(q)_1}\\
  &=\frac{\mathbf{H}^2}{(q)_1}+\mathbf{H} \, \mathbf{G}_r\Big(2+\frac{2q-1}{(q)_1}-\frac{1}{(q)_1}\Big)\\
  &=\frac{\mathbf{H}^2}{(q)_1}.\\
  \end{align*}
  Thus for any integer $r\geq 2$ the cardinal $F_{r}(n)$ of $\mathcal{F}_{r}(n)$ is equal to the coefficient of $q^n$ in $\mathbf{H}^2/(q)_1$. This means that $F_r(n)=F_{r+1}(n)$ for any integer $r\geq 2$.  
   \end{proof}
   \begin{rem} As we have seen in the proof of Theorem \ref{main}, we don't really need to compute $S_{3,b}$. However, we have done it in Lemma \ref{S3}.
   \end{rem}
   
 
\bibliographystyle{acm}
\bibliography{Afsharijoo_Gonzalez_Mourtada}

\hspace{1cm}

\address{Instituto de Matem\'atica Interdisciplinar, Departamento de \'Algebra, Geometr\'ia y Topolog\'ia,
Facultad de Ciencias Matem\'aticas, Universidad Complutense de Madrid, Plaza de las Ciencias 3,
Madrid 28040, Espa\~na.\\
\email{pafshari@cyu.fr} \\
\email{pgonzalez@mat.ucm.es}\\}

\address{Universit\'e Paris Cit\'e, Sorbonne Universit\'e, CNRS, Institut de Math\'ematiques de Jussieu-Paris Rive Gauche, Paris, 75013, France,\\
\email{hussein.mourtada@imj-prg.fr} }

\end{document}